\documentclass[12pt]{amsart}
\usepackage{amssymb, amsmath, amsthm}

\newtheorem{lemma}{Lemma}[section]
\newtheorem{prop}[lemma]{Proposition}
\newtheorem{cor}[lemma]{Corollary}
\newtheorem{thm}[lemma]{Theorem}

\theoremstyle{definition}
\newtheorem{rem}[lemma]{Remark}

\newtheorem{defi}[lemma]{Definition}
\newtheorem{ex}[lemma]{Example}

\newcommand{\slth}{SL_3({\mathbb H})}
\newcommand{\slto}{E_{6(-26)}}
\newcommand{\sltc}{SL_3({\mathbb C})}

\newcommand{\sltr}{SL_3({\mathbb R})}

\newcommand{\cat}{\text{CAT}}
\newcommand{\Set}[1]{\mathcal{#1}}
\newcommand{\nbd}[2]{\mathcal{N}_{#2} ({#1}) }
\newcommand{\bb}{\mathbb}

\begin{document}
\title{Some Geometric Groups with Rapid Decay}
\author{Indira Chatterji$^\dagger$ and Kim Ruane$^{\dagger\dagger}$}
\address{Math Department, Cornell University, Ithaca NY, USA and Math Department, Tufts University, Somerville MA, USA}
\email{indira@math.cornell.edu, kim.ruane@tufts.edu}
\date{\today}
\thanks{$^\dagger$ Partially supported by the Swiss National Funds}
\thanks{$^{\dagger\dagger}$ Partially supported by Tufts University FRAC award}
\begin{abstract}We explain simple methods to establish the property of Rapid Decay for a number of groups arising geometrically. Those lead to new examples of groups with the property of Rapid Decay, notably including non-cocompact lattices in rank one Lie groups.\end{abstract}
\maketitle
\section*{Introduction}
A discrete group $\Gamma$ is said to have {\it the property of Rapid Decay (property RD) with respect to a length function $\ell$} if there exists a polynomial
$P$ such that for any $r\in{\bb R}_+$ and any $f$ in the complex group algebra ${\bb C}\Gamma$ supported on elements of length shorter than $r$ the following inequality holds:
$$\|f\|_*\leq P(r)\|f\|_2$$
where $\|f\|_*$ denotes the operator norm of $f$ acting by left convolution on $\ell^2(\Gamma)$, and $\|f\|_2$ is the usual $\ell^2$ norm. Property RD had a first striking application in A. Connes and H. Moscovici's work proving the Novikov conjecture for Gromov hyperbolic groups \cite{CM} and is now relevant in the context of the Baum-Connes conjecture, mainly due to the work of V. Lafforgue in \cite{Lafforgue2} (see Section \ref{appl}). First established for free groups by U. Haagerup in \cite{Haagerup}, property RD has been introduced and studied as such by P. Jolissaint in \cite{Jolissaint}, who notably established it for groups of polynomial growth, and for classical hyperbolic groups. The extension to Gromov hyperbolic groups is due to P. de la Harpe in \cite{Harpe}. The first examples of higher rank groups with property RD have been given by J. Ramagge, G. Robertson and T. Steger in \cite{RRS}, where they established it for $\tilde{A}_2$ and $\tilde{A}_1\times\tilde{A}_1$ groups. V. Lafforgue proved property RD for cocompact lattices in $\sltr$ and $\sltc$ in \cite{Lafforgue}. His result has been generalized by the first author in \cite{moi} to cocompact lattices in $\slth$ and $\slto$ as well as in a finite product of rank one Lie groups. It is well-known (see Section \ref{basics}) that non-cocompact lattices in higher rank simple Lie groups do not have property RD, and it is a conjecture due to Valette that all cocompact lattices in a semisimple Lie group should have property RD (see~\cite{Valette}). In this paper we shall see that the situation is different in rank one.  Indeed, all lattices have property RD. More precisely we prove the following.
\begin{thm}\label{cool}Groups which are hyperbolic relative to a family of polynomial growth subgroups satisfy property RD.\end{thm}
This result has recently been generalized in \cite{DS2}. The following is then an immediate consequence.
\begin{cor}\label{Rk1}
\begin{itemize}
\item[(a)]Let $M$ be a complete and simply connected Riemannian manifold of pinched negative curvature. Any discrete and finite covolume subgroup of Isom($M$) has property RD. In particular, all lattices in rank one Lie groups have property RD.
\item[(b)]Suppose $G$ acts properly discontinuously, cocompactly, by isometries on a CAT(0) space with the Isolated Flats Property. Then $G$ has property RD.
\end{itemize}
\end{cor}
Due to the work of Lafforgue in \cite{Lafforgue2}, the following is then straightforward.
\begin{cor}\label{bc}\begin{itemize}
\item[(a)]Let $M$ be a complete and simply connected Riemannian manifold of pinched negative curvature and bounded curvature tensor. Any discrete and finite covolume subgroup of Isom($M$) satisfies the Baum-Connes conjecture. In particular, all lattices in rank one Lie groups satisfy the Baum-Connes conjecture.
\item[(b)]Suppose $G$ acts properly discontinuously, cocompactly, by isometries on a CAT(0) space with the Isolated Flats Property.  Then $G$ satisfies the Baum-Connes Conjecture. \end{itemize}\end{cor}
Closed subgroups (and in particular lattices) in $SO(n,1)$ and $SU(n,1)$ were known to satisfy the Baum-Connes conjecture due to the work of Julg and Kasparov in \cite{JK} on the Baum-Connes conjecture with coefficients. The case of cocompact lattices in rank one Lie groups follows from the work of Lafforgue in \cite{Lafforgue2} (see Skandalis' exposition \cite{SB}) and was a break-through in the subject because it provided the first examples of property (T) groups satisfying the Baum-Connes conjecture. Closed subgroups of $Sp(n,1)$ and of the exceptional Lie group $F_{4(-20)}$ are due to a recent work of Julg in \cite{J} and \cite{J2} where he proves the Baum-Connes conjecture with coefficients for those groups. General word hyperbolic groups were shown to satisfy Baum-Connes in \cite{YuMin}. The role of property RD in Lafforgue's work will briefly be recalled in Section \ref{appl}. 
Combining the computations used to prove Theorem \ref{cool} with further geometrical considerations, we also prove the following.
\begin{thm}\label{RDcube}Groups acting properly with uniformly bounded stabilizers and cellularly on a CAT(0) cube complex of finite dimension have property RD.\end{thm}
The particular case where the CAT(0) cube complex is an arbitrary finite product of trees was treated independently in \cite{moi} and by M. Talbi in \cite{TalbiCR} and \cite{Talbi}. The latter also discusses the case (for which property RD is still open in general) of groups acting cocompactly on euclidean buildings.  Here, he obtains interesting geometric information leading to partial results. The Baum-Connes Conjecture was already known for groups acting on a CAT(0) cube complex using the work of \cite{NibloReeves2} combined with the work of Higson and Kasparov \cite{HK} on a-(T)-menable groups satisfying the Baum-Connes conjecture (see also the exposition of P. Julg in \cite{JB} and \cite{welches} for a-(T)-menable groups). However, we mention the following consequence which follows immediately from Jolissaint's work (see Corollary 3.1.8 in \cite{Jolissaint} or Theorem \ref{obstruction} below).
\begin{cor}\label{NoAmSup}Groups acting properly with uniformly bounded stabilizers on a CAT(0) cube complex of finite dimension cannot contain amenable subgroups of super-polynomial growth.\end{cor}
Recently, D. Wise and M. Sageev in \cite{WS} proved a version of the Tits alternative which is stronger than the above. 

The paper is organized as follows. In Section \ref{basics} we recall some basics regarding property RD and give a crucial geometric condition which is sufficient to imply property RD (Proposition \ref{max}). Section \ref{relhyp} is devoted to relatively hyperbolic groups and the proof of Theorem \ref{cool}. We establish property RD by showing that the groups of Theorem \ref{cool} satisfy the assumptions of Proposition \ref{max}. In Section \ref{appl} we discuss Corollary \ref{Rk1}, as well as the Baum-Connes conjecture and the applications of our results in this context, explaining Corollary \ref{bc}. Finally, Section \ref{cubes} deals with CAT(0) cube complexes and the proof of Theorem \ref{RDcube}. Again we establish property RD by showing that the groups of Theorem \ref{RDcube} satisfy the assumptions of Proposition \ref{max}.
\section{Rapid Decay and techniques}\label{basics}
We will explain the basic notions related to property RD for discrete groups. Except for Proposition \ref{max}, the results given in this section are either simple remarks or results contained in Jolissaint's paper \cite{Jolissaint}. There is also a theory of property RD for locally compact groups, first studied in \cite{Jolissaint} and developed further in \cite{JiSch} and \cite{RDloc}, but we restrict here to discrete groups.
\begin{defi}
Let ${G}$ be a group, a \emph{length function} on ${G}$ is a map $\ell:{G}\to{\bb R}_+$ satisfying: 
\begin{itemize}
\item $\ell(e)=0$, where $e$ denotes the neutral element in ${G}$,
\item $\ell(\gamma)=\ell(\gamma^{-1})$ for any $\gamma\in{G}$,
\item $\ell(\gamma\mu)\leq \ell(\gamma)+\ell(\mu)$ for any $\gamma,\mu\in{G}$.\end{itemize}\end{defi}
If ${G}$ is generated by some finite subset $S$, then the \emph{algebraic word length} $L_S:{G}\to{\bb N}$ is a length function on ${G}$, where, for $\gamma\in{G}$, $L_S(\gamma)$ is the minimal length of $\gamma$ as a word on the alphabet $S\cup S^{-1}$, namely,
$$L_S(\gamma)=\min\{n\in{\bb N}|\gamma=s_1\dots s_n,\,s_i\in S\cup S^{-1}\}.$$
For a length function $\ell$, the map $d_{\ell}(\gamma,\mu)=\ell(\gamma^{-1}\mu)$ is a left ${G}$-invariant pseudo-distance on ${G}$. We will write $B_\ell(\gamma,r)$ for the ball of center $\gamma\in{G}$ and radius $r$ with respect to the pseudo-distance $d_{\ell}$, and simply $B(\gamma,r)$ when the context is clear.
\begin{defi} Denote by ${\bb C}{G}$ the complex group algebra of $G$, that we view as the set of functions $f:{G}\to{\bb C}$ with finite support, denoted supp($f$). The ring structure is given by pointwise addition and convolution:
$$f*g(\gamma)=\sum_{\mu\in{G}}f(\mu)g(\mu^{-1}\gamma),$$
(for $f,g\in{\bb C}{G}$ and $\gamma\in{G}$). We denote by ${\bb R}_+{G}$ the subset of ${\bb C}{G}$ consisting of functions with target in ${\bb R}_+$. We shall consider the following completions of ${\bb C}{G}$:
\begin{itemize}
\item[(a)] the {\it reduced C*-algebra of ${G}$}, given by $C^*_r(G)=\overline{{\bb C}{G}}^{\|\ \|_*}$, where $\|f\|_{*}=\sup\{\|f*g\|_2\,|\,\|g\|_2=1\}$ is the \emph{operator norm} of $f\in{\bb C}{G}$,
\item[(b)] for $s\geq 0$, the $s$-Sobolev space $H^s_\ell({G})=\overline{{\bb C}{G}}^{\|\ \|_{\ell,s}}$, where $\|f\|_{\ell,s}=\sqrt{\sum_{\gamma\in{G}}|f(\gamma)|^2(1+\ell(\gamma))^{2s}}$ is a {\it weighted $\ell^2$ norm}.  For $s=0$, this is $\ell^2({G})$, the Hilbert space of square summable functions on ${G}$.
\end{itemize}
\end{defi}
\begin{defi}[P. Jolissaint, \cite{Jolissaint}]\label{deBase} Let $\ell$ be a length function on ${G}$. We say that ${G}$ has the {\it property of Rapid Decay}\footnote{Some authors refer to those groups as satisfying the \emph{Haagerup inequality}. That sometimes leads non-experts to confusion with the \emph{Haagerup property} or \emph{a-(T)-menability}, something very different, see \cite{welches}.} with respect to $\ell$, if there exists $C,s>0$ such that, for each $f\in{\bb C}{G}$ one has
$$\|f\|_{*}\leq C\|f\|_{\ell,s}.$$
\end{defi}
The functions in the intersection of all Sobolev spaces
$$H^{\infty}_\ell({G})=\bigcap_{s\geq 0}H^s_\ell({G})$$
are called \emph{rapidly decaying} functions, as their decay at infinity is faster than any inverse of a polynomial in $\ell$. Property RD with respect to $\ell$ is equivalent to having $H^{\infty}_\ell({G})\subseteq C^*_r({G})$ (see Remark 1.2.2 in \cite{Jolissaint} or Section 2 of \cite{RDloc} for a more detailed proof), which explains the terminology. In case where ${G}={\bb Z}$, one checks that under Fourier transform $C^*_r({\bb Z})$ is isomorphic to the C*-algebra of continuous functions over the circle, and $H^{\infty}_\ell({G})$ corresponds to smooth functions. The following will give us more flexibility in the computations.
\begin{prop}\label{RD}Let ${G}$ be a discrete group, endowed with a length function $\ell$. Then the following are equivalent:
\begin{itemize}
\item[(1)] The group ${G}$ has property RD with respect to $\ell$.
\item[(2)] There exists a polynomial $P$ such that, for any $r>0$ and any $f\in{\bb R}_+{G}$ so that supp($f$) is contained in $B_{\ell}(e,r)$
$$\|f\|_*\leq P(r)\|f\|_2.$$
\item[(3)] There exists a polynomial $P$ such that, for any $r>0$ and any $f,g,h\in{\bb R}_+{G}$ so that supp($f$) is contained in $B_{\ell}(e,r)$
$$f*g*h(e)\leq P(r)\|f\|_2\|g\|_2\|h\|_2.$$
\item[(4)] Any subgroup $H$ in ${G}$ has property RD with respect to the induced length.
\end{itemize}\end{prop}
\begin{proof} We start with the equivalence between (1) and (2). Take $f\in{\bb C}{G}$ with support contained in $B_{\ell}(e,r)$, we have:
$$\|f\|_*\leq C\|f\|_{\ell,s}\leq C\sqrt{\sum_{\gamma\in B(e,r)}|f(\gamma)|^2(r+1)^{2s}}=C(r+1)^s\|f\|_2$$
and thus (2) is satisfied, for the polynomial $P(r)=C(r+1)^s$. Conversely, if we denote by $S_n=\{\gamma\in{G}|n\leq\ell(\gamma)< n+1\}$, for $n\in{\bb N}$ we compute, for $f\in{\bb R}_+{G}$:
\begin{eqnarray*}\|f\|_*&=&\|\sum_{n=0}^\infty f|_{S_n}\|_*\leq\sum_{n=0}^\infty\|f|_{S_n}\|_*\leq\sum_{n=0}^\infty P(n+1)\|f|_{S_n}\|_2\\
&\leq&\sum_{n=0}^\infty C(n+1)^k\|f|_{S_n}\|_2=C\sum_{n=0}^\infty (n+1)^{-1}(n+1)^{k+1}\|f|_{S_n}\|_2\\
&\leq&C\sqrt{\sum_{n=0}^\infty(n+1)^{-2}}\sqrt{\sum_{n=0}^\infty (n+1)^{2k+2}\|f|_{S_n}\|_2^2}
=C\frac{\pi}{\sqrt{6}}\|f\|_{\ell,k+1}.\end{eqnarray*}
We finish by noticing that for $f\in{\bb C}{G}$, if one writes $f=f_1-f_2+i(f_3-f_4)$ with $f_i\in{\bb R}_+{G}$ and the supports of $f_i$ and $f_{i+1}$ are disjoint for $i=1,3$, then $\|f\|^2_2=\sum_{i=1}^4\|f_i\|^2_2$ and thus
$$\|f\|_*\leq\sum_{i=1}^4\|f_i\|_*\leq P(r)\sum_{i=1}^4\|f_i\|_2\leq P(r)\sqrt{4}\|f\|_2.$$

\smallskip

Let us turn to the equivalence between (2) and (3). To see that (3) implies (2) it is enough to define, for $\gamma\in{G}$
$$h(\gamma)=\frac{f*g(\gamma^{-1})}{\|f*g\|_2},$$
and notice that, in that case, $f*g*h(e)=\|f*g\|_2$ and $\|h\|_2=1$. We deduce (2) decomposing $g=g_1-g_2+i(g_3-g_4)$ as above. That (2) implies (3) follows from the Cauchy-Schwartz inequality:
$$f*g*h(e)=\sum_{\gamma\in{G}}(f*g)(\gamma)\check{h}(\gamma)\leq\|f*g\|_2\|h\|_2\leq P(r)\|f\|_2\|g\|_2\|h\|_2$$
where $\check{h}(\gamma)=h(\gamma^{-1})$.

\smallskip

Finally, (4) implies (1) trivially since ${G}$ is a subgroup of itself (and the induced length is the original one) and (2) implies (4) since if $H$ is a subgroup of ${G}$, $f\in{\bb R}_+H$ supported in a ball of radius $r$ can be viewed in ${\bb R}_+{G}$, supported in a ball of radius $r$ as well, thus
$$\|f\|_{*,H}\leq\|f\|_{*,{G}}\leq P(r)\|f\|_{\ell^2{G}}=P(r)\|f\|_{\ell^2H}.$$
\end{proof}
Recall that a discrete group ${G}$ has {\it polynomial growth with respect to a length $\ell$} if there exists a polynomial $P$ such that the cardinality of the ball of radius $r$ (denoted by $|B_{\ell}(e,r)|$) is bounded by $P(r)$. Combined with point (4) of the previous proposition, the following result gives the only known obstruction to property RD, namely the presence of an amenable subgroup of superpolynomial growth.
\begin{thm}[P. Jolissaint, Corollary 3.1.8 in \cite{Jolissaint}]\label{obstruction} Let ${G}$ be a discrete amenable group and $\ell$ a length function on ${G}$, the following are equivalent. 
\begin{itemize}
\item[(i)] ${G}$ has property RD with respect to $\ell$.
\item[(ii)] ${G}$ is of polynomial growth with respect to $\ell$. 
\end{itemize}
Moreover, the growth will be bounded by $P^2$, if $P$ is the polynomial of Proposition \ref{RD} (2).\end{thm}
\begin{proof}[Proof (taken from \cite{Valette})] We use a variation of Kesten's characterization of amenability (see \cite{L}), stating that a group ${G}$ is amenable if and only if $\|f\|_1=\|f\|_*$ for any $f\in{\bb R}_+{G}$.

$(i)\Rightarrow (ii)$: Let $f$ be the characteristic function of $B_{\ell}(e,r)$ and assume property RD. We use Proposition \ref{RD} (2):
$$P(r)\|f\|_2\geq\|f\|_*=\|f\|_1=|B_{\ell}(e,r)|=\sqrt{|B_{\ell}(e,r)|}\|f\|_2.$$
$(ii)\Rightarrow (i)$: Take $f\in{\bb C}{G}$ such that supp$(f)\subset B_{\ell}(e,r)$, then:
$$\|f\|_*\leq\|f\|_1=\sum_{\gamma\in B_{\ell}(e,r)}|f(\gamma)|\leq\sqrt{|B_{\ell}(e,r)|}\sqrt{\sum_{\gamma\in B_{\ell}(e,r)}|f(\gamma)|^2}\leq P(r)\|f\|_2,$$
the second inequality being just the Cauchy-Schwartz inequality, and the last inequality is polynomial growth.\end{proof}
According to A. Lubotzky, S. Mozes and M. S. Raghunathan in \cite{LMR} there exists an infinite cyclic subgroup growing exponentially with respect to the word length in any non-cocompact lattice in higher rank (exponentially distorted copy of ${\bb Z}$), and hence Theorem \ref{obstruction} combined with Proposition \ref{RD} (4) shows that non-cocompact lattices in higher rank Lie groups cannot have property RD. An important point of the present paper is to show that it is not the case for non-cocompact lattices in rank one Lie groups. It is part of a conjecture due to A. Valette (see \cite{Valette}, Conjecture 7) that cocompact lattices in semisimple Lie groups should have property RD.
\begin{rem} It is well-known and easy (see Lemma 1.1.4 in \cite{Jolissaint}) that a finitely generated group ${G}$ has property RD with respect to the word length as soon as it has property RD for any other length, and thus explains why we will be sloppy regarding the length functions involved as soon as we deal with finitely generated groups.
\end{rem}
 The following proposition is a reformulation of Proposition 2.3 in \cite{Lafforgue} and it will be our main tool to prove property RD in this paper.
\begin{prop}\label{max} Let ${G}$ be a group acting freely and by isometries on a metric space $(X,d)$ such that there is a ${G}$-equivariant map
\begin{eqnarray*}C: X\times X &\to & {\mathcal P}(X)\\
(x,y)&\mapsto& C(x,y)\end{eqnarray*}
(where ${\mathcal P}(X)$ are the subsets of $X$) satisfying the following (for any $x,y,z\in X$ and $\gamma\in{G}$):
\begin{itemize}
\item[(i)] $C(x,y)\cap C(y,z)\cap C(z,x)\not=\emptyset$.
\item[(ii)] There is a polynomial $P$ such that for any $r\geq 0$, then the cardinality of $C(x,y)\cap B(x,r)$ is bounded above by $P(r)$.
\item[(iii)] There is a polynomial $Q$ such that if $d(x,y)\leq r$, then the diameter of $C(x,y)$ is bounded by $Q(r)$.\end{itemize}
Then ${G}$ has property RD (with respect to the length $\ell(\gamma)=d(x,\gamma x)$, $x\in X$ any base point).\end{prop}
\begin{proof}
Let us consider the groupoid ${\mathcal G}$ given as follows:
$${\mathcal G}=X\times X/\sim$$
where $(x,y)\sim (s,t)$ if there exists $\gamma\in{G}$ with $x=\gamma s,y=\gamma t$. We write $[x,y]$ for the class of the pair $(x,y)$ in ${\mathcal G}$, and
$${\mathcal G}^0=\{[x,y]\in{\mathcal G}| x=\gamma y\ \hbox{for some }\gamma\in{G}\}$$
with source and range given by
\begin{eqnarray*}{s,r:{\mathcal G}}& \to & {\mathcal G}^0\\
{[x,y]} & \mapsto & s[x,y]=[y,y],r[x,y]=[x,x]\end{eqnarray*}
so that the composable elements are
$${\mathcal G}^2=\{([x,y],[s,t])\in{\mathcal G}\times{\mathcal G}|y=\gamma s\ \hbox{for a }\gamma\in{G}\}$$
and $[x,y]\cdot [s,t]=[x,y]\cdot [\gamma s,\gamma t]=[x,\gamma t]$ if $y=\gamma s$. For $f,g\in{\bb R}_+{\mathcal G}$, the convolution is given by
$$f*_{\mathcal G}g[x,y]=\sum_{z\in X}f[x,z]g[z,y],$$
(for $[x,y]\in{\mathcal G}$). It is enough to prove that there exists a polynomial $P$ such that for every $r\in{\bb R}_+$ and every $f,g,h\in{\bb R}_+{\mathcal G}$ such that ${\rm supp}(f)\subset{\mathcal G}_r=\{[x,y]\in{\mathcal G}|d(x,y)\leq r\}$, then the following inequality holds:
\begin{equation}\label{RDgroupoide}f*_{\mathcal G}g*_{\mathcal G}h[x,x]\leq P(r)\|f\|_2\|g\|_2\|h\|_2\end{equation}
for every $x\in X$, where $\|f\|_2^2=\sum_{[x,y]\in{\mathcal G}}f[x,y]^2$. Indeed, from (1) we conclude that ${G}$ has property RD by using Proposition \ref{RD} (4) and defining for a fixed $x_0\in X$, a linear map $T:{\bb C}{G}\to{\bb C}{\mathcal G}$ by 
$$T(f)[x,y]=\left\{\begin{array}{cc}f(\gamma)&\hbox{ if }[x,y]=[x_0,\gamma x_0]\\
0&\hbox{ otherwise },\end{array}\right.$$
so that $T(f)[x_0,x_0]=f(e)$. One checks that $\|T(f)\|_2=\|f\|_2$ and that $T(f*g)=T(f)*_{\mathcal G}T(g)$ for any $f,g\in{\bb C}{G}$, and hence $T$ is an isometric embedding of algebras. We now turn to the proof of inequality (\ref{RDgroupoide}) above. For $x_0\in X$, 
\begin{eqnarray*}& &f*_{\mathcal G}g*_{\mathcal G}h[x_0,x_0]=\sum_{y,z\in X^2}f[x_0,y]g[y,z]h[z,x_0]\\
&\leq&\sum_{x\in{G}\setminus X}\sum_{y,z\in X^2}f[x,y]g[y,z]h[z,x]=\sum_{x,y,z\in{G}\setminus X^3}f[x,y]g[y,z]h[z,x]\end{eqnarray*} 
and because of assumption (i), we have that
$$\sum_{x,y,z\in{G}\setminus X^3}f[x,y]g[y,z]h[z,x]\leq \sum_{\stackrel{x,y,z,t\in{G}\setminus X^4}{t\in C(x,y)\cap C(y,z)\cap C(z,x)}}f[x,y]g[y,z]h[z,x].$$
For $t\in C(x,y)$ and $d(x,y)\leq r$, then $d(x,t)\leq Q(r)$ and $d(t,y)\leq Q(r)$ by assumption (iii). Let $H_1=\ell^2{\mathcal G}$, $H_2=\ell^2{\mathcal G}_{Q(r)}=H_3$ and denote by $\{[t,x]|[t,x]\in{\mathcal G}\}$ (respectively $\{[t,x]|[t,x]\in{\mathcal G}_{Q(r)}\}$) the canonical basis for $H_1$ (respectively $H_2$ and $H_3$). Now we define $T_f\in{\mathcal L}(H_2,H_3)$ given on the basis vectors of $H_2$ and $H_3$ by
$$\left<T_{f}([t,x]),[v,y]\right>_{H_3}=\left\{\begin{array}{cc}f[x,y]&\hbox{ if }t\hbox{ is in the orbit of }v\hbox{ (so we assume $t=v$),}\\
&\hbox{ and if }t\in C(x,y)\\
0&\hbox{ otherwise.}\end{array}\right.$$
We then extend $T_f$ by linearity to an element of ${\mathcal L}(H_2,H_3)$. In the same way we define $T_g\in{\mathcal L}(H_1,H_3)$ and $T_h\in{\mathcal L}(H_2,H_1)$. For $[t,x]\in{\mathcal G}_{Q(r)}$ we have
$$\left<T_f\circ T_g\circ T_h([t,x]),[t,x]\right>_{H_2}=\sum_{\stackrel{y,z\in X^2}{t\in C(x,y)\cap C(y,z)\cap C(z,x)}}f[x,y]g[y,z]h[z,x]$$
(note that this equality uses condition (iii) and the fact that $f$ is supported on a ball of radius $r$), thus
\begin{eqnarray*}\hbox{Trace}(T_f\circ T_g\circ T_h)&=&\sum_{[x,t]\in{\mathcal G}_{Q(r)}}\left<T_f\circ T_g\circ T_h([t,x]),[t,x]\right>_{H_2}\\
&=&\sum_{\stackrel{x,t\in{G}\setminus X^2}{d(x,t)\leq Q(r)}}\sum_{\stackrel{y,z\in X^2}{t\in C(x,y)\cap C(y,z)\cap C(z,x)}}f[x,y]g[y,z]h[z,x]\\
&=&\sum_{\stackrel{x,y,z,t\in{G}\setminus X^4}{t\in C(x,y)\cap C(y,z)\cap C(z,x)}}f[x,y]g[y,z]h[z,x]\end{eqnarray*}
(the last equality again uses condition (iii) and the fact that $f$ is supported on a ball of radius $r$). Now we use that $\hbox{Trace}(T_f\circ T_g\circ T_h)\leq \|T_f\|_{HS}\|T_g\|_{HS}\|T_h\|_{HS}$ (this holds for Hilbert-Schmidt operators in general) and evaluate those Hilbert-Schmidt norms:
\begin{eqnarray*}\|T_f\|^2_{HS}&=&\sum_{([t,x],[v,y])\in {\mathcal G}_{Q(r)}^2}|\left<T_{f}([t,x]),[v,y]\right>_{H_3}|^2=\sum_{\stackrel{x,y,t\in{G}\setminus X^3,d(x,t)\leq Q(r)}{t\in C(x,y)}}|f[x,y]|^2\\
&\leq& P(Q(r))\sum_{x,y\in{G}\setminus X^2}|f[x,y]|^2=P(Q(r))\|f\|^2_2,\end{eqnarray*}
the last inequality holding because of assumption (ii). Similarly, one shows that $\|T_g\|^2_{HS}\leq P(Q(r))\|g\|^2_2$ and $\|T_h\|^2_{HS}\leq P(Q(r))\|h\|^2_2$.\end{proof}
\begin{rem}If a group $G$ has polynomial growth, one checks that taking $X=G$ and defining $C(x,y)=B(x,d)\cup B(y,d)$ where $d=d(x,y)$ is the word length associated to any generating set for $G$, fulfills the assumptions of Proposition \ref{max} above.  This reproves the fact, due to P. Jolissaint in \cite{Jolissaint}, that polynomial growth groups have property RD.

If a group $G$ is Gromov hyperbolic, then any Cayley graph is $\delta$-hyperbolic for some $\delta\geq 0$, meaning that for any geodesic triangle with vertices $x,y,z$ in the Cayley graph, the geodesic between two of the vertices is contained in the $\delta$-neighborhood of the union of the two other geodesics (see \cite{GdlH}, Chapitre 1, D\'efinition 27). Let us take $X=G$ with a distance induced by the length associated to a finite generating set $S$ and define $C(x,y)$ as to be the set of elements $g\in G$ that are in the $\delta$-neighborhood of all geodesics between $x$ and $y$, where $\delta\geq 0$ is the hyperbolicity constant for the length associated to the generating set $S$. One then checks that the assumptions of Proposition \ref{max} above are fulfilled, reproving the fact, due to P. de la Harpe and P. Jolissaint in \cite{Harpe}, that Gromov hyperbolic groups have property RD.

A combination of those two cases will give a proof for Theorem \ref{cool}.\end{rem}
\section{Relatively hyperbolic groups}\label{relhyp}
This section will be devoted to the proof of Theorem ~\ref{cool}. Let us first notice that a lattice $G$ in a 
rank one Lie group acts properly discontinuously by 
isometries on a rank one symmetric space $X$.  The space $X$ is negatively curved 
as a Riemannian manifold, but also $\delta$-hyperbolic in the 
sense of Gromov (indeed, the classical notion of curvature is stronger 
than Gromov's $\delta$-hyperbolicity).  If $G$ is cocompact, then it is a Gromov hyperbolic group and has property RD, but a non-cocompact lattice $G$ will not  necessarily be Gromov hyperbolic as the parabolic subgroups will sometimes be free abelian. However, $G$ will be hyperbolic relative to the family of  {\it parabolic subgroups}.   

\smallskip  

The notion of a group $G$ being hyperbolic relative to a collection of subgroups $\{H_1,\dots , H_k\}$ was first introduced by Gromov in \cite{Gr}, Section 8.6.  Recall that a Gromov hyperbolic group can be  defined as a group which acts properly discontinuously by isometries and cocompactly on a proper (closed metric balls are compact), geodesic,  $\delta$-hyperbolic metric space.  One could loosely define a group $G$ to be hyperbolic relative to a collection of subgroups $\{H_1,\dots , H_k\}$ by dropping the cocompactness assumption and replacing it by the requirement that the quotient be quasi-isometric to the union of $k$ copies of rays $[0,\infty)$ joined at 0. The subgroups $\{H_1,\dots , H_k\}$ are the isotropy groups of the end of each ray (see details below). This is the view point taken in Gromov's original definition. 

\smallskip  

One can also define a group $G$ to be Gromov hyperbolic if one (and hence any) Cayley graph of $G$ is a $\delta$-hyperbolic metric space. In \cite{F}, Farb gives an alternative definition of relative hyperbolicity in terms of properties of a {\it coned-off} Cayley graph (see below for formal definitions). The Gromov and Farb definitions are not equivalent as is shown in \cite{Sz}.  However, according to Bowditch's work in \cite{Bo}, Farb's definition together with the Bounded Coset Penetration Property (see Definition~\ref{BCP} below) is equivalent to the Gromov definition. More on relatively hyperbolic groups can be found in \cite{Dahmani} or \cite{DS1}.

\smallskip  

We shall use both the Gromov and Farb definitions and start here by describing more precisely Gromov's viewpoint. Let $X$ be a proper, geodesic, $\delta$-hyperbolic metric space, recall (from Section III.3 of \cite{BH}) that one defines $\partial{X}$, the boundary of $X$, as equivalence classes of geodesic rays in $X$. For $c:[0,\infty)\to X$ a geodesic ray, the {\it Busemann function of $c$} is the map
\begin{eqnarray*}b_c:X\to{\bb R}, x\mapsto b_c(x):=\limsup_{t\to\infty} d(x,c(t))-t.\end{eqnarray*}
According to Remark 3.4 in Section III.3 of \cite{BH} one can construct $\partial{X}$, the boundary of $X$, as equivalence classes of Busemann functions.
\begin{defi}\label{horofunction}For $c$ a geodesic ray and $r\geq 0$, the sublevel sets $b_c^{-1}(-\infty,r]\subset X$ are called (closed) {\it horoballs} of radius $r$ and the sets $b_c^{-1}(r)$ are called {\it horospheres centered at} $\xi=c(\infty)\in\partial{X}$.  
\end{defi}

Now suppose that a group $G$ acts properly discontinuously by isometries on $X$ so that the quotient of $X$ by $G$ is quasi-isometric to the union of $k$ copies of $[0,\infty)$ joined at zero.  For simplicity, assume the action is free. Lift the rays in the quotient to obtain $k$ points $p_1,p_2,\ldots ,p_k$ in $\partial X$.  Choose geodesic rays $c_i:[0,\infty)\to X$ such that $c_i(\infty)=p_i$ for $i=1,2,\ldots ,k$.  Let $H_i$ be the isotropy subgroup for the action of $G$ on $\partial X$ of $c_i(\infty)=p_i$ and assume that $H_i$ preserves the Busemann function $b_{c_i}$. If there exists a $G$-invariant system of horoballs ${\mathcal B}$ centered at the points $p_1,\ldots ,p_k$, in $X$ so that the action of $G$ on $X\setminus(\cup_{B\in{\mathcal B}} B)$ is cocompact, then we say the action of $G$ on $X$ is {\it geometrically finite with parabolic subgroups} $\{H_1,\dots,H_k\}$. Note that, since each $H_i$ stabilizes a distinct point of the boundary, the intersection $H_i\cap H_j$ is finite if $i\not= j$, and hence up to taking finite index subgroups we can assume that  $\{H_1,\dots,H_k\}$ pairwise intersect trivially.
\begin{defi}[Gromov]\label{Gromov} A group $G$ is called \emph{hyperbolic relative to the family $\{H_1,\dots,H_k\}$}
of finitely generated subgroups if $G$ admits a geometrically finite action on a proper, geodesic $\delta$-hyperbolic metric space $X$ with parabolic subgroups $\{H_1,\dots,H_k\}$.
\end{defi}
\begin{ex}To understand the Gromov definition, one should think of a non-cocompact lattice $G$ in $SO(n,1)$,  the isometry group of $\mathbb H^n$. In this case, the group $G$ acts properly discontinuously by isometries on $\mathbb H^n$ but not cocompactly - the quotient is a manifold (or orbifold if the action is not free) with finitely many cusps.  However, $G$ does act cocompactly on a subset of $\mathbb H^n$.  Indeed, one can remove a $G$-invariant set of disjoint open horoballs about the parabolic fixed points of $G$ in $\partial\mathbb H^n$ and $G$ acts cocompactly on the complement of this in $\mathbb H^n$.   It is this viewpoint that the Gromov definition generalizes.   \end{ex}
Recall that  a subset $Z$ of $X$ is called {\it quasi-convex} (or {\it $K$-quasi-convex})  if there exists a constant $K>0$ such that for any points $x,y\in Z$, any geodesic in $X$ from $x$ to $y$ is in the $K$-neighborhood of $Z$. In the case of a $\delta$-hyperbolic space $X$, the horoballs are quasi-convex subsets of $X$ and thus the horospheres are quasi-convex in the complement of the union of the horoballs.

\begin{defi}For $r\geq 0$, a system of quasi-convex subsets of a metric space $X$ is called {\it $r$-separated} if $d(Q_1, Q_2)>r$ for any pair of sets $Q_1,Q_2$ in the collection.\end{defi}  

The following result, due to Bowditch, says that  if $G$ is hyperbolic relative to the family $\{H_1,\dots,H_k\}$, then one can make a careful choice of horoballs so that the action of each $H_i$ on the corresponding horosphere is cocompact.
\begin{lemma}[Lemma 6.3 and 6.12 in \cite{Bo}]\label{bowditch}  If $p_i\in\partial X$ is the parabolic point stabilized by $H_i$, then there is a $G$-invariant, $r$-separated system of horoballs $\mathcal B$ such that $G$ acts cocompactly on $X'=X\setminus\cup_{B\in\mathcal B}B$.   Furthermore, suppose that $B_i\in\mathcal B$ is a horoball stabilized by the subgroup $H_i$, then $H_i$ acts cocompactly on the bounding horosphere $S_i$ in $X$.  
\end{lemma}
The analogous lemma in the case of $\mathbb H^n$ says that one can equivariantly shrink the horoballs in $\mathbb H^n$ so that they are sufficiently far apart to guarantee the action of a parabolic group on the corresponding horosphere is cocompact. In the classical case, the horoballs are convex subsets of $\mathbb H^n$ and thus the horospheres are convex in the complement of the union of the horoballs.

\medskip

For Farb's definition, we begin with a finitely generated group $G$ with a fixed generating set $S$ and a finite set of infinite, finitely generated subgroups $H_1,\dots ,H_k$ of $G$ that pairwise intersect trivially.   Consider  the Cayley graph $\Gamma=\Gamma(G,S)$ of $G$ with respect to $S$ and the usual right action of $G$ on $\Gamma$ by multiplication.  Add a vertex $c_{gH_i},\ i=1,2,\ldots ,k$ for each left coset $gH_i$ of $H_i$ in $G$, and connect $c_{gH_i}$ with each $x\in gH_i$ by an edge of length $1\over 2$.  The new graph is denoted by $\hat\Gamma=\hat\Gamma(H_1,H_2\ldots ,H_k)$ and is called the {\it coned-off} Cayley graph of $G$ with respect to $\{H_1,H_2,\ldots ,H_k\}$.  We denote by $\hat d$ the path metric on $\hat\Gamma$.   It is easy to see that $\hat\Gamma$ is quasi-isometric to the graph obtained from $\Gamma$ by identifying each left coset to a point. 
\begin{defi}\label{weakhyp} We call the group $G$ \emph{weakly hyperbolic relative to $\{H_1,H_2,\ldots ,H_k\}$} if $\hat\Gamma$ is a $\delta$-hyperbolic metric space.  
\end{defi}
This is Farb's original definition for $G$ being hyperbolic relative to the collection $\{H_1\ldots ,H_k\}$, however here we use the terminology suggested by Bowditch in \cite{Bo}.   We will now describe the Bounded Coset Penetration Property (see also Section 3.3 in \cite{F}).

To simplify matters, we restrict for a while to the case of having only one parabolic subgroup $H$ in $G$. This is similar to considering a hyperbolic manifold with one cusp versus several.  The statements and proofs are easier to read in this setting and we explain at the end of the section how to handle more than one subgroup in the collection. We also assume that the generating set $S$ for $G$ contains a generating set $S_H$ for $H$.  For a word $z$ in the letters of $S$, denote by $\overline z$ the group element obtained as the endpoint of the path in $\Gamma$ whose initial point is the identity of $G$ and follows the edge labels given in the word $z$. A path $w$ in $\Gamma$ is a word $z$ in the letters of $S$, together with an initial point $x_0$, so that the endpoint of the path is the element $\overline{z}x_0\in G$. Given a path $w$ in $\Gamma$, we obtain a path $\hat w$ in $\hat\Gamma$ as follows:  read $w$ from right to left and identify maximal subwords $z$ in the generators $S_H$.  If $z$ is a maximal $S_H$ subword from a vertex $g$ to $g\overline z$ in $\Gamma$, we can replace the subpath $z$ by an edge path with two edges each of length ${1\over 2}$ in $\hat\Gamma$ - namely,  one edge from the vertex $g$ to the cone point $c_{gH}$ and an another edge from $c_{gH}$ to the vertex $g\overline z$.  Following Farb, these subpaths $z$ of a word $w$ are called {\it coset subwords}.  The correspondence $w\mapsto\hat w$ is clearly a surjective map.  If $\hat w$ passes through some cone point $c_{gH}$, then we say  that $w$ {\it penetrates} the coset $gH$.  

Recall that for $P\geq 1$, a {\it $P$-quasi-geodesic} in a geodesic metric space $(X,d)$ is a $P$-quasi-isometric embedding $\varphi:[a,b]\to X$, i.e. $\varphi$ is a map satisfying
$$\frac{1}{P}|t-t'|\leq d(\varphi(t),\varphi(t'))\leq P|t-t'|$$
(where $a<t\leq t'<b$). We call the image of $\varphi$ a $P$-quasi-geodesic between $\varphi(a)$ and $\varphi(b)$. If $\hat{w}$ is a geodesic (or $P$-quasi-geodesic) in $\hat{\Gamma}$, then we call any preimage $w$ of $\hat{w}$ a {\it relative geodesic} (or {\it relative $P$-quasi-geodesic}).  A path $w$ in $\Gamma$ (or $\hat w$ in $\hat\Gamma$) is called a {\it path without backtracking} if, for every coset $gH$ which $\hat w$ penetrates, $\hat w$ never returns to that coset after leaving. 
\begin{defi}[B. Farb]\label{BCP} The pair $(G,H)$ is said to satisfy the {\it Bounded Coset Penetration Property} (BCP for short) if,  for every $P\ge 1$, there is a constant $K=K(P)\geq 0$ so that if $u$ and $v$ are relative $P$-quasi-geodesics in $\hat\Gamma$ without backtracking that start at the same point and end in the same coset, then the following are true:
\begin{itemize}
\item[1.]If $u$ penetrates a coset and $v$ does not penetrate that coset, then $u$ traveled a $\Gamma$-distance of at most $K$  in that coset. 
\item[2.] If they both penetrate a coset, then the $\Gamma$-distance between their entry and
exit points is at most $K$ (but they can travel a long time in that coset).
\end{itemize} 
\end{defi}
To generalize Farb's definition to the case of $G$ hyperbolic relative to a finite family of subgroups $\{H_1,\dots,H_k\}$, we assume that the subgroups pairwise intersect trivially and we take a generating set $S$ for $G$ that contains $S_{H_1}\cup\dots\cup S_{H_k}$, where $S_{H_i}$ is a generating set for $H_i$. Since the $S_{H_i}$'s are pairwise disjoint, decomposing a word in coset subwords makes sense and Definition \ref{BCP} has a straightforward generalization. We now have both definitions of relative hyperbolicity. According to Theorem 7.10 of \cite{Bo}, Gromov's definition (here Definition \ref{Gromov}) is equivalent to Farb's definition (here Definition ~\ref{weakhyp}) with the Bounded Coset Penetration Property. 

\smallskip

The following shows how the subgroups $\{H_1,\dots,H_k\}$ sit in a relatively hyperbolic group $G$.  
\begin{lemma}\label{QI} For $G$ hyperbolic relative to $\{H_1,\dots,H_k\}$, there is $M\geq 0$ such that the image of the inclusion map from $H_i$ (with its induced word metric) to $G$ is $M$-quasi-convex (for any $i=1,\dots,k$).  As a consequence, each $H_i$ is quasi-isometrically embedded in $G$.
\end{lemma}
\begin{proof} The idea is to think of Gromov's Definition \ref{Gromov} first and use Lemma ~\ref{bowditch} to identify $H_i$ with a horosphere $S_i$ in $X'=X\setminus\cup_{B\in\mathcal B}B$, as well as the fact that $S_i$ is quasi-convex inside of $X'$.  Then, since the Cayley graph $\Gamma=\Gamma(G,S)$ is quasi-isometric to $X'$ we can carry the information about $H_i$ back to $\Gamma$. More precisely, the Cayley graph $\Gamma$ sits inside $X'$ via $f:\Gamma\to X'$, $g\mapsto g\cdot x_0$, where $x_0\in S_i\subset X'$ is a fixed base-point. For $s\in S$, the element $s\cdot x_0$ lies in $s\cdot S_i$ which is a horosphere disjoint from $S_i$ if $s$ is not in $S_{H_i}$.   If $s\in S_{Hi}$, then $s\cdot x_0$ lies in $S_i$.  Extend the map equivariantly and we see that each coset $gH_i$ lies in its own horosphere $g\cdot S_i$. The map $f$ is a quasi-isometry.  

The first observation is that there is a constant $L\geq 0$ such that for any $h_1,h_2\in H_i$, we can find an $L$-quasi-geodesic path $\beta$ in $\Gamma$ from $h_1$ to $h_2$ that is contained in the subgroup $H_i$. Indeed, in the space $X'$, the points $x_1=h_1\cdot x_0$ and $x_2=h_2\cdot x_0$ are both in $S_i$ which is quasi-convex in $X'$.  Hence  there exists a constant $L\geq 0$ such that any $X'$-geodesic $\gamma$ from $x_1$ to $x_2$ lies in the $L$-neighborhood of $S_i$.  We also know from Lemma ~\ref{bowditch} that $H_i$ acts cocompactly on $S_i$.  Thus there exists $N\geq 0$ such that every point of $S_i$ is within $N$ of an $H_i$-orbit point.  Now move along $\gamma$ at unit speed, and at each integer $t$ between $0$ and $d(h_1\cdot x_0,h_2\cdot x_0)$, we can pick an $H_i$-orbit point that is $(N+L)$-close to $\gamma(t)$. This sequence of $H_i$-orbit points will give a path in $\Gamma$ from $h_1$ to $h_2$ using only vertices from $H_i$. Let us call $\beta$ this path, which is a quasi-geodesic because $f$ is a  quasi-isometry  between $\Gamma$ and $X'$.  

We want to show that there exists a constant $M\geq 0$ such that any $\Gamma$-geodesic $\alpha=[h_1,h_2]$ between elements $h_1,h_2$ of $H_i$ is in the $M$-neighborhood of the subgroup $H_i$. Suppose that $\alpha$ does not penetrate the coset $H_i$.  Then by the BCP, Definition \ref{BCP} point 1, the length of $\beta$ would be at most $K=K(L)$.  But the length of $\alpha$ is certainly less than the length of $\beta$ since $\alpha$ is geodesic.  Thus $\alpha$ has length no more than $K$.  This means every point of $\alpha$ is within $K\over 2$ of $H_i$.  If $\alpha$ does penetrate $H_i$, then we decompose $\alpha$ into successive pieces that do and do not penetrate $H_i$.  Each piece lying outside $H_i$ has length at most $K\over 2$ by the above argument.  Thus $M={K\over 2}$ is the desired constant.  
\end{proof}

We proceed with proving that a group $G$ which is hyperbolic relative to $\{H_1,\dots,H_k\}$ satisfies the hypotheses of Proposition ~\ref{max}, where the metric space $X$ is just $G$ with the Cayley graph metric. We start with the construction of a map 
$$C:G\times G\to {\mathcal P}(G)$$ 
that satisfies the conditions of Proposition~\ref{max}. To simplify matters, we again restrict for a while to the case of having only one parabolic subgroup $H$ in $G$. The following is the first step.
\begin{defi}Let $G$ be a finitely generated group and $H$ a finitely generated subgroup of $G$. For $\delta\geq 0$ and $x,y$ in $G$, we define 
$${V}_{\delta}(x,y)=\{t\in G|d(x,t)+d(t,y)\leq d(x,y)+\delta\}\subseteq G,$$
and $\hat{V}_{\delta}(x,y)=\{c_{gH}|\hbox{ such that }{V}_{\delta}(x,y)\cap gH\ne\emptyset\}$.
\end{defi}
In other words, ${V}_{\delta}(x,y)$ is a $\delta$-thickening of the convex hull of $x,y$ (i.e. of the set of all geodesics between $x$ and $y$), and  $\hat{V}_{\delta}(x,y)$ consists of the cone points of $\hat{\Gamma}$  which are in the $\delta$-neighborhood of a $\Gamma$-geodesic from $x$ to $y$. The crucial property of those sets $\hat{V}_{\delta}(x,y)$ in our context is the following.
\begin{lemma}\label{AlmostL} If $G$ is hyperbolic relative to a subgroup $H$, then there is a $\delta$ big enough so that for any $x,y,z\in G$, then 
$$\hat{V}_{\delta}(x,y)\cap\hat{V}_{\delta}(y,z)\cap\hat{V}_{\delta}(z,x){\ne} \emptyset.$$\end{lemma}
\begin{proof}According to Theorem 1.12 (3) of \cite{DS1}, quasi-geodesics in the Cayley graph of $G$ are at bounded Hausdorff distance from any geodesic in the coned-off graph $\hat{\Gamma}$. The lemma now follows from the fact that the coned-off graph $\hat{\Gamma}$ is hyperbolic.\end{proof}
\begin{defi}Let $G$ be a finitely generated group and $H$ be a finitely generated subgroup of $G$. For $r\geq 0$ and $x$ in the coset $xH$, we denote by $B_H(x,r)=B(x,r)\cap xH$ (this is the intersection of the $xH$ points of the ball of radius $r$ centered at $x$). For $x,y\in G$ we define 
$$C_H(x,y)=\left\{\begin{array}{cc}B_H(x,r)\cup B_H(y,r) & \hbox{ if } xH=yH\hbox{ with }r=d(x,y)+K\\
\{x\}\cup\{y\} &\hbox{ otherwise.}\end{array}\right.$$ 
and $K$ is a constant (that we will later choose to be the BCP constant of Definition \ref{BCP}).\end{defi} 
Obviously for three points $x,y,z$ in the same $H$-coset, then $C_H(x,y)\cap C_H(y,z)\cap C_H(z,x)\not=\emptyset$. Keeping Lemma \ref{QI} in mind, the crucial property of the sets $C_H(x,y)$ in our context is the following.
\begin{lemma}\label{NotFarFromH} Let $G$ be a finitely generated group and $H$ be polynomial growth subgroup of $G$ that is quasi-isometrically embedded. For any $x,y\in G$, then the cardinality of $B(x,r)\cap C_H(x,y)$ is bounded by the growth polynomial of $H$.\end{lemma}
\begin{proof} In case where $x$ and $y$ are not in the same $H$-coset then there is nothing to prove. Otherwise, observe that $B(x,r)\cap C_H(x,y)\subseteq B_H(x,r)=x(B(1,r)\cap H)$.  Since $H$ has polynomial growth and the growth polynomial is a quasi-isometry invariant in the class of discrete groups the lemma now follows.\end{proof}
Now we can define the sets $C(x,y)$ for $x$ and $y$ arbitrary elements of $G$: Take any quasi-geodesic $\gamma=\gamma(x,y)$ in $\Gamma$ between $x$ and $y$ such that all its vertices belong to ${V}_{\delta}(x,y)$ (we then say that $\gamma\subset {V}_{\delta}(x,y)$). Looking at maximal $H$-coset subwords as explained just before Definition \ref{BCP}, we can determine the coset penetration points of $\gamma$. We write them as follows (this includes entrance and exit points):
$$x=x_0^{\gamma},x_1^{\gamma},\dots,x_n^{\gamma}=y.$$
We finally define
$$C(x,y)=\bigcup_{\gamma\subset {V}_{\delta}(x,y)}\bigcup_{1\leq i,j\leq n}C_H(x_i^{\gamma},x_j^{\gamma})\subseteq G.$$
These sets are $G$-invariant by construction.
\begin{proof}[Proof of Theorem \ref{cool}] We shall prove that the above defined map 
$$C:G\times G\to {\mathcal P}(G)$$
satisfies the conditions of Proposition~\ref{max}. Let us start by proving point (i). Suppose $x,y,z\in G$ and assume that $d(x,y)\geq\max\{d(y,z),d(z,x)\}$. If $x,y,z$ belong to the same coset, then $z\in C(x,y)$ since $d(x,z)\leq d(y,x)$ and hence $z\in C(x,y)\cap C(y,z)\cap C(z,x)\not=\emptyset$. If they all belong to different cosets, then according to Lemma \ref{AlmostL}, there are 3 relative $\delta$-quasi-geodesics $\gamma(x,y),\gamma(y,z)$ and $\gamma(z,x)$ in $\Gamma$ entering a common coset $gH$. Assume that $\gamma(x,y),\gamma(y,z)$ and $\gamma(z,x)$ enter and leave the coset $gH$ at respective points $a_{xy},b_{xy}$, $a_{yz},b_{yz}$ and $a_{zx},b_{zx}$.  Note, that the points $b_{xy}$ and $a_{yz}$ (respectively $b_{yz}$ and $a_{zx}$ as well as $b_{zx}$ and $a_{xy}$) are at most distance $K$ apart by Definition \ref{BCP} point 2., hence the intersection
$$C_H(a_{xy},b_{xy})\cap C_H(a_{yz},b_{yz})\cap C_H(a_{zx},b_{zx})\not=\emptyset$$
and is contained in the intersection of $C(x,y),C(y,z)$ and $C(z,x)$, showing that the latter is nonempty as well. If only two out of the three points belong to the same coset, then we can assume that $xH=yH\not=zH$. There are quasi-geodesics $\gamma(x,z)$ and $\gamma(y,z)$ in $\Gamma$ entering the coset $xH=yH$ at respective points $a$ and $b$ which lie at a bounded distance $K$ by Definition \ref{BCP} point 2.  If $d(x,y)\geq\max\{d(a,x),d(b,y)\}$, then both $a$ and $b$ are in $C_H(x,y)$, hence $C(x,y)\cap C(y,z)\cap C(z,x)\not=\emptyset$. If $d(a,x)\geq\max\{d(x,y),d(b,y)\}$, then $y\in C_H(x,a)$ so that again $C(x,y)\cap C(y,z)\cap C(z,x)\not=\emptyset$.

\smallbreak

We now turn to point (ii) and take $x,y\in G$. To start with, since $G$ is finitely generated $|\hat{V}_{\delta}(x,y)\cap B(x,r)|\leq Cr$, where $C$ is the number of cosets points in $\hat{\Gamma}$ contained in a ball of radius $\delta$ centered at a coset in $\hat{\Gamma}$.  We know that $H$ has polynomial growth (both with its own word length or with the length induced by $G$ as $H$ is quasi-isometrically embedded in $G$, see Lemma~\ref{QI}) and that according to Lemma \ref{NotFarFromH} all the intersections $C_H(x_i^{\gamma},x_j^{\gamma})\cap B(x,r)$ have cardinality bounded by a polynomial in $r$.  The BCP property implies there are at most $B(e,K)$ of the sets $C_H(x_i^{\gamma},x_j^{\gamma})\cap B(x,r)$, which proves that $C(x,y)$ has cardinality bounded by a polynomial as well.  

\smallbreak

For point (iii), it is enough to show that there is a constant $C$ such that for any $x,y$ at distance less than $r$ and for any $z\in C(x,y)$ then $d(z,x)\leq Cr+\delta$. Take $z\in C(x,y)$, then either $z\in{V}_{\delta}(x,y)$, or $z$ is in a $C_H(a,b)$ for some $a,b$ on $\gamma\subset{V}_{\delta}(x,y)$. For $z\in\hat{V}_{\delta}(x,y)$, then ${d}(z,x)\leq r+\delta$ by definition. For a $z\in C_H(a,b)$ we have $d(z,x)\leq d(z,a)+d(a,x)\leq 2r+\delta$ because for any $a,b\in\gamma$, then $d(a,b)\leq r+\delta$.

\medbreak

For the general case, few things need to be changed. For each subgroup $H_i$, one can define the sets $C_i(x,y)$ as done above. Then one defines 
$$C(x,y)=\bigcup_{i=1}^nC_i(x,y).$$
That condition (i) is satisfied follows from from the fact that the coned off graph is hyperbolic, and hence if the 3 points $x,y,z$ lie in different cosets, there are 3 relative $\delta$-quasi-geodesics entering a common coset $gH_i$ (analogue to Lemma \ref{AlmostL}) and the discussion is the same as above otherwise. That conditions (ii) and (iii) hold follows from the fact that they hold in each coned off graph (relatively to an $H_i$), independently of relative hyperbolicity but only depend on the BCP.
\end{proof}
\begin{rem}Intuitively, for $x,y\in G$ the set $C(x,y)$ is constructed as follows. One fixes $\delta$ big enough and takes a $G$-orbit $G\cdot x_0$ in the hyperbolic space. One then considers a geodesic between $x\cdot x_0$ and $y\cdot x_0$. When this geodesic intersects a horoball, one modifies it by replacing the segment in the horoball by a ball in the horosphere, with diameter the distance between the entry and exit points. The set $C(x,y)$ consists of all elements in $G$ whose image in the $G$-orbit $G\cdot x_0$ is within distance $\delta$ of a modified geodesic between $x\cdot x_0$ and $y\cdot x_0$.\end{rem}
\section{Applications}\label{appl}
In this section, we explain Corollaries \ref{Rk1} and \ref{bc} stated in the introduction.   
\begin{proof}[Proof of Corollary \ref{Rk1}(a)] Fundamental groups of finite volume quotients of Riemannian manifolds with pinched negative sectional curvature are known to have finitely many cusps (see \cite{Eb} Lemma 3.1d), and the cusps are quasi-isometric to a ray. This, according to Gromov's Definition \ref{Gromov}, shows that such groups are hyperbolic relative to their (finitely many) cusp subgroups.   According to \cite{Eb}, Corollary 3.3, the cusp subgroups are virtually nilpotent and thus have polynomial growth. This includes discrete subgroups of Isom($M$) with finite volume quotient, where $M$ is a non-compact, simply connected, real rank one symmetric space.  These are exactly the lattices in rank one Lie groups.
\end{proof}
The proof of Corollary ~\ref{Rk1}(b) only require the relevant definitions.   The following notion was first evident in the work of \cite{KapovichLeeb95} on 3-manifolds and was expanded upon in C. Hruska's thesis, \cite{Hruska2ComplexIFP} in the setting of CAT(0) 2-complexes.  Recently, Hruska and Kleiner in \cite{HL}  extended the work in \cite{Hruska2ComplexIFP} to all higher dimensions.  We recall that a \emph{flat} in a metric space is an isometrically embedded copy of ${\bb R}^n$.
\begin{defi}[Isolated Flats Property]\label{def:IFP}
A CAT(0) metric space $X$ has the \emph{Isolated Flats Property} (or IFP for short) if it contains a family of flats~$\Set{F}$ with the following two properties.
\begin{enumerate}
\item There is a constant $C$ so that every flat in $X$ lies in the $C$-neighborhood of some flat $F \in \Set{F}$.
\item There is a function $\psi \colon \mathbb R_{+} \to \mathbb R_{+}$ such that for any two distinct flats $F_1, F_2 \in \Set{F}$ and for any positive number $r$, the intersection
$$\nbd{F_1}{r} \cap \nbd{F_2}{r}$$
of Hausdorff neighborhoods of $F_1$ and $F_2$ has diameter at most $\psi(r)$.
\end{enumerate}
\end{defi}
If we consider two maximal flats to be equivalent when their Hausdorff
distance is finite, then the family $\Set{F}$
in the preceding definition consists of one maximal flat from each
equivalence class.
Intuitively, $X$ has the Isolated Flats Property if given any two
maximal flats in $X$ which are not parallel, the two flats
diverge from each other in all directions.
In particular, two maximal flats are either parallel, or
\emph{disjoint at infinity}, meaning that their corresponding
boundary spheres are disjoint.

\begin{proof}[Proof of Corollary ~\ref{Rk1}(b)] According to C. Hruska and B. Kleiner \cite{HL}, groups acting properly by isometries and cocompactly on CAT(0) spaces with the Isolated Flats Property are hyperbolic relative to the stabilizers of the flats.  These subgroups are virtually abelian hence of polynomial growth so that Theorem ~\ref{cool} applies.  In the particular where $X$ is a CAT(0) 2-complex with IFP, that follows from earlier work of Hruska in \cite{Hruska2ComplexIFP}.\end{proof}
Corollary \ref {bc} is straightforward from Lafforgue's work \cite{Lafforgue2} and Corollary ~\ref{Rk1}. In order to explain that, we recall a few facts relating property RD and the Baum-Connes conjecture. We will not attempt to describe what the Baum-Connes conjecture is but we just state it. To do so, we denote by $\underline{E}G$ the classifying space for proper actions of a discrete group $G$ and by $K_{i}^{G}(\underline{E}G)$ its equivariant $K$-homology (for $i=0,1$). P. Baum and A. Connes in \cite{BaumConnes} defined an \emph{assembly map} $\mu_i:K_{i}^{\Gamma}(\underline{E}G)\to K_{i}(C^{*}_{r}(G))$ ($i=0,1$), where $K_{i}(C^{*}_{r}(G))$ is the topological $K$-theory of the C*-algebra $C^{*}_{r}(G)$, and formulated the following.

\medbreak
\noindent
{\bf Conjecture} (Baum-Connes). {\it Let $\Gamma$ be a discrete group. The assembly map
$$\mu_{i}^{\Gamma}:K_{i}^{G}(\underline{E}G)\rightarrow K_{i}(C^{*}_{r}(G))$$
is an isomorphism.}

\medbreak
We refer to \cite{Valette} or \cite{MV} for an introduction to the Baum-Connes conjecture, and to \cite{JB} or \cite{SB} for expositions of Higson-Kasparov or Lafforgue crucial advances on this conjecture. The first evidence of the relevance of property RD in the context of the Baum-Connes conjecture was given by the work of Connes and Moscovici \cite{CM} on the very closely related Novikov conjecture that we won't discuss here. The following observation is due to Lafforgue, Proposition 1.2 of \cite{Lafforgue}, generalizing a theorem by A. Connes, explained by Jolissaint in \cite{JolK}, Theorem A.
\begin{prop}[V. Lafforgue]\label{obs} Assume that $G$ has property RD with respect to $\ell$. Then $H^s_\ell(G)$ is a Banach algebra for $s$ large enough and the inclusion $H^s_\ell(G)\hookrightarrow C^*_r(G)$ induces an isomorphism in $K$-theory.
\end{prop}
\begin{defi}[V. Lafforgue] A Banach algebra ${\mathcal A}{G}$ is an \emph{unconditional completion} of ${\bb C}{G}$ if it contains ${\bb C}{G}$ as a dense subalgebra and if, for $f_1,f_2\in{\bb C}{G}$ such that $|f_1(\gamma)|\leq|f_2(\gamma)|$ for all $\gamma\in{G}$, we have
$$\|f_1\|_{{\mathcal A}{G}}\leq\|f_2\|_{{\mathcal A}{G}}.$$\end{defi}
The reduced C*-algebra is in general not an unconditional completion, even for ${G}={\bb Z}$. If ${G}$ has property RD with respect to a length function $\ell$, then for $s$ large enough, $H^s_\ell({G})$ is a convolution algebra and an unconditional completion. For any unconditional completion ${\mathcal A}{G}$ of ${\bb C}{G}$, Lafforgue constructs a map 
$$\mu_{\mathcal A}:K_{*}^{{G}}(\underline{E}{G})\to K_{*}({\mathcal A}{G}),$$ 
compatible with the assembly map $\mu_*$. He then defines a class ${\mathcal C}'$ of groups, closed by products, containing (among others) discrete groups acting properly discontinuously and isometrically either on a simply connected Riemannian manifold with non positive curvature bounded from below and whose curvature tensor has bounded derivatives, or properly discontinuously, isometrically and cocompactly on a CAT(0) metric space (more precisely, properly discontinuously and isometrically on a weakly geodesic and strongly bolic metric space, see \cite{Lafforgue2} for the definition of strong bolicity) and proves the following.
\begin{thm}[V. Lafforgue \cite{Lafforgue2}] For any group belonging to the class ${\mathcal C}'$ and for any unconditional completion ${\mathcal A}{G}$ of ${\bb C}{G}$ the map $\mu_{\mathcal A}$ is an isomorphism.
\end{thm}
In view of Proposition \ref{obs} above, this gives the following.
\begin{thm}[V. Lafforgue \cite{Lafforgue2}]\label{C+RD=BC} The Baum-Connes conjecture holds for any property RD group belonging to the class ${\mathcal C}'$.\end{thm}
We refer to Skandalis' nice exposition in \cite{SB} for more on Lafforgue's contribution to the Baum-Connes conjecture. Since the groups of Corollary \ref{bc} are contained in Lafforgue's class ${\mathcal C}'$, this finishes the proof of Corollary \ref{bc}. We end this section by giving several examples of CAT(0) spaces with the Isolated Flats Property.
\begin{ex}  Suppose that $M$ is a Haken 3-manifold obtained by gluing hyperbolic components along torus boundary components.  Then according to Kapovich and Leeb in \cite{KapovichLeeb95}, the universal cover $\widetilde{M}$ of $M$ has the Isolated Flats Property.
\end{ex}  

\begin{ex}  Any finite covolume discrete subgroup of Isom($\mathbb H^n$) acts geometrically on a CAT(0) space with the Isolated Flats Property (see Proposition 9.1 in \cite{Hruska2}).\end{ex}

\begin{ex}
D. Wise has shown that a proper, cocompact piecewise Euclidean
CAT(0) 2-complex $X$ has the Isolated Flats Property if and only
if $X$ does not contain an isometrically embedded triplane.  In particular, this includes any $\cat(0)$ 2-complex built out of hexagons.   A \emph{triplane} is a space formed by isometrically gluing three Euclidean half planes together along their boundary lines. For a proof, see \cite{Hruska2ComplexIFP}.  
\end{ex}
\begin{rem} It has been shown in \cite{YuMin} that Gromov hyperbolic groups belong to Lafforgue's class ${\mathcal C'}$. However, it is an an open question which relatively hyperbolic groups can be given a metric which is strongly bolic, and hence the Baum-Connes conjecture for the groups of Theorem \ref{cool} is still open.\end{rem}
\section{CAT(0) cube complexes}\label{cubes}
In this section we give the proofs of Theorem ~\ref{RDcube} and Corollary ~\ref{NoAmSup}. We shall prove that the zero-skeleton of a cubical CAT(0) complex, endowed with the distance of the one-skeleton, satisfies the hypothesis of Proposition \ref{max}. We start by recalling some definitions.
\begin{defi}Let $(X,d)$ be a metric space and $\delta\geq 0$. For any finite sequence of points $x_1,\dots,x_n\in X$, we say that $x_1\dots x_n$ is a {\it $\delta$-path} if
$$d(x_1,x_2)+\dots+d(x_{n-1},x_n)\leq d(x_1,x_n)+\delta$$
and that three points $x,y,z\in X$ form a {\it $\delta$-retractable triple} if there exists $t\in X$ such that the paths $xty$, $ytz$ and $ztx$ are $\delta$-paths. In this case we will say that the triple $x,y,z$ $\delta$-retracts on $t$, and that $t$ is a \emph{$\delta$-midpoint} or an \emph{approximate midpoint}. We will say that $X$ satisfies {\it property $L_\delta$}, or is an \emph{$L_{\delta}$-space} if there exists a $\delta\geq 0$ such that any triple has a $\delta$-midpoint. Notice that if a triple is $\delta$-retractable, then it is $\delta'$-retractable for any $\delta'\geq\delta$.\end{defi}
\begin{rem}Hyperbolic metric spaces are $L_{\delta}$-spaces, and ``to be an $L_{\delta}$-space'' is closed under direct product (with an $\ell^1$-combination of the distances) but not under quasi-isometries. Some considerations of these spaces can be found in \cite{these} and they will be the object of an independent study in \cite{Ldelta}. Groups admitting a Cayley graph which is an $L_{\delta}$-space have sub-cubic isoperimetric inequality as shown by Elder in \cite{Elder}.\end{rem}
\begin{defi}A \emph{cube complex} $X$ is a metric polyhedral complex in which each cell is isometric to the Euclidean cube $[-1/2,1/2]^n$ and the gluing maps are isometries. A cube complex $X$ is called CAT(0) if the metric induced by the Euclidean metric on the cubes gives $X$ a CAT(0) metric. We shall denote by $X^i$ the set of $i$-dimensional cells of $X$, and say that $X$ is finite dimensional if there is $n<\infty$ such that $X^m$ is empty for any $m>n$.\end{defi}
We shall use the following fundamental work developed by Sageev in \cite{Sag}, that we now recall. A \emph{combinatorial hyperplane} is an equivalence class of unoriented edges, where two edges $e$ and $f$ are called equivalent if there exists a finite sequence of edges $e=e_1,\dots,e_n=f$, such that for each $i=1,\dots,n-1$, $e_i$ and $e_{i+1}$ are opposite sides of some 2-cube in $X$. A \emph{hyperplane} is a subcomplex formed by cells (=subcubes) in the first barycentric subdivision and which are orthogonal (when viewed in the non-subdivided cube as a subset of ${\bb R}^n$) to the edges of a combinatorial hyperplane. The crucial result we shall need is the following.
\begin{thm}[Sageev, \cite{Sag} Theorem 4.6]\label{Sag} Given two vertices $p$ and $q$ in $X^0$ that are at distance $n$ in the 1-skeleton, any geodesic path crosses $n$ distinct hyperplanes. Moreover, each of those hyperplanes separates $p$ and $q$, and any edge path between $p$ and $q$ must cross these hyperplanes.\end{thm}
In other words the distance between two points in the 1-skeleton only depends on the number of hyperplanes separating them. This already yields the following useful fact.
\begin{lemma}\label{lacetPair}Any closed loop in the 1-skeleton has even length.\end{lemma}
\begin{proof}Follows from Theorem \ref{Sag} above, since in a closed loop, every hyperplane has to be crossed an even number of times.\end{proof}
\begin{defi}For a given closed loop in the 1-skeleton we define its \emph{combinatorial area} by the minimal number of two cells (squares) needed to fill the loop. This number exists because we assumed $X$ to be CAT(0), hence contractible.\end{defi}
We now proceed with a result that, as pointed out by G. Niblo, was already known by M. Roller in \cite{South} in the context of median graphs.  We give a proof here for completeness.  
\begin{prop}\label{CATestL}The zero-skeleton $X^0$ of a CAT(0) cube complex $X$, endowed with the distance of the one-skeleton, is an $L_0$-space.\end{prop}
\begin{proof}The first step is to show that a triple $x,y,z\in X^0$ with $d(x,y)=2$ and $d(x,z)=d(y,z)=n$ retracts.\\ There are two cases to be considered, the first one is the case where $x$ and $y$ belong to a common two-dimensional cube, case in which they have to be opposite vertices. We call $a$ and $b$ the two remaining vertices (opposite as well in this common two-dimensional cube). By Lemma \ref{lacetPair}, $d(a,z)$ and $d(b,z)$ cannot be $n$ -- e.g. $a,z,y$ lie on a closed loop -- so it is either $n+1$ or $n-1$ since both $a$ and $b$ are adjacent to $x$ and $y$. If $d(a,z)=n-1$ we are done, $xaz$, $yaz$ and $xay$ are 0-paths, so let us assume that $d(a,z)=n+1$. This means that the hyperplane spanned by the equivalence class of the oriented edge from $x$ to $a$ separates $a$ from $z$ but not $b$ from $z$. Similarly for the hyperplane spanned by the equivalence class of the oriented edge from $y$ to $a$, so that between $b$ and $z$ there are two hyperplanes less than between $a$ and $z$, and we deduce that $d(b,z)\leq d(a,z)-2=n-1$, hence $xbz$, $ybz$ and $xby$ are 0-paths.\\ The second case to consider is where $d(x,y)=2$ and $x$ and $y$ do not share a common 2-dimensional cube, so that there is a unique element $t$ at distance one from both $x$ and $y$. Again because of Lemma \ref{lacetPair}, the distance between $t$ and $z$ is either $n+1$ or $n-1$. If it's $n-1$ we are done, and it cannot be $n+1$ because then no geodesic from $x$ to $z$ and from $y$ to $z$ would go through $t$, which means that $x,t,y,z,x$ would form a non contractible closed loop.\\ 
To finish the proof we proceed by contradiction and assume that there is a triple that doesn't retract. Among all those triples let $x,y,z$ be one with smallest possible minimum of the three side lengths, say $d(x,y)=n$ ($n$ is bigger than $2$ by the preceding discussion), and consider the geodesics between $x,y$ and $z$ realizing the smallest possible combinatorial area (and keep them for the rest of the proof). On the geodesic between $x$ and $y$ pick $a$ at a distance one from $x$, so that $d(a,y)=n-1$ and hence the triple $a,y,z$ retracts, on a point $t$. It is easy to see that actually $t=y$: indeed, if $t\not=y$, then $x,t,z$ would give a non-retractable triple with strictly smaller minimum side length than the triple $x,y,z$ (indeed, if the triple $x,t,z$ was to retract, then so would the triple $x,y,z$). Since the triple $a,y,z$ retracts on $y$, it means that the path $ayz$ is a geodesic. We know that $d(a,z)=d(x,z)\pm 1$ (because $a$ is at distance $1$ to $x$ and closed loops have even length), and it cannot be that $d(a,z)=d(x,z)-1$, because then $a$ would lie on a geodesic between $x$ and $z$, and hence the triple $x,y,z$ would retract on $y$ as well. Hence  $d(a,z)=d(x,z)+1$, and we now are almost reduced to the first step of the proof: Take $b$ on the chosen geodesic from $x$ to $y$ at distance $1$ to $y$, by assumption on the minimality of the triple $x,y,z$, the triple $x,b,z$ has to retract, on a point $t$ which is easily seen to be at distance $1$ from $b$. We get a contradiction because now the triple $t,y,z$ falls in the first step of the proof ($t=y$ is not possible because we assumed $b$ on a geodesic to $x$), hence is contractible, which allows to contract the triple $x,y,z$ as well.\end{proof}
In order to proceed with the proof of Theorem \ref{RDcube} we will need the following result of \cite{Sag}, which says that a collection of pairwise intersecting hyperplanes have to share a common cube:
\begin{thm}[Sageev, \cite{Sag} Theorem 4.14]\label{intersection} For $X$ a CAT(0) cube complex, if $h_1,\dots,h_k$ is a collection of hyperplanes such that $h_i\cap h_j\not=0$, then $\bigcap_{i=1}^kh_i\not=0$.\end{thm}
We now proceed with the proof of Theorem \ref{RDcube}.
\begin{proof}[Proof of Theorem \ref{RDcube}]We will first show that $X^0$ satisfies the hypothesis of Proposition \ref{max}, and this will settle the case where $G$ acts freely on the CAT(0) cube complex. We will treat the general case at the end of the proof. The map $C:X^0\times X^0\to {\mathcal P}(X^0)$ is simply defined as follows:
$$C(x,y)=\{t\in X^0\hbox{ such that }d(x,t)+d(t,y)=d(x,y)\}.$$
This map is $G$-equivariant follows from the $G$-invariance of the distance $d$. That this map satisfies (i) follows from Lemma \ref{CATestL}, point (iii) is obvious (since $C(x,y)$ consists of points on geodesics) so let us prove point (ii). We have to show that for any $x,y$ in $X^0$, the number of points $t$ in $X^0$ on a geodesic from $x$ to $y$ and lying in a ball of radius $r$ centered at $x$ is at most polynomial, and the polynomial is actually of degree $n$, the dimension of the cube complex. The idea is to show that the largest number of points occurs when $x$ and $y$ are opposite vertices of an $n$-cube in $\mathbb R^n$ with side length $r\over 2$ .  This number is clearly bounded by a polynomial of degree $n$ in $r$.

Let $H$ be the set of hyperplanes separating $x$ from $y$. We call two hyperplanes $h$ and $h'$ parallel if they don't intersect, and write $H$ as a disjoint union of subsets $P_1,\dots,P_k$ where all hyperplanes in a given $P_i$ are parallel. A partition is called minimal if for each $i$ and each $h\in P_i$ there exists $j$ and $h'\in P_j$ intersecting $h$.  On one extreme, if a minimal partition has just one piece, then there is a unique geodesic between $x$ and $y$.  On the other extreme, if there are $n$ pieces each containing one wall, then $x$ and $y$ are opposite vertices in an $n$-cube in $\bb R^n$.

We now claim that a minimal partition has at most $n$ pieces: Take a geodesic from $x$ to $y$ and define $P_1$ as follows; put in $P_1$ the first hyperplane crossed, say $h_1$. Put $h_2$ in $P_1$ if and only if it doesn't intersect $h_1$, so that at the $i$-th step one puts $h_i$ in $P_1$ if and only if $h_i$ intersects none of the hyperplanes already in $P_1$. Then define $P_i$ similarly, starting with the first hyperplane not already in $P_{i-1}$ and skipping all the hyperplanes already sitting in $P_{i-1}$. Doing so there is a sequence $\{h_i\}$ with $h_i\in P_i$ of pairwise intersecting hyperplanes and hence according to Sageev (Theorem \ref{intersection} cited above), this partition (which is not unique as it depends on the geodesic we started with) has at most $n$ pieces.

Let us now treat the general case and produce a metric space $Y$ on which $G$ acts freely and satisfying the assumptions of Proposition \ref{max}. The metric space is 
$$Y=X^0\coprod_{x\in X^0}G_x,$$
where $G_x<G$ is the stabilizer of $x\in X^0$, and the metric is given as follows: two distinct points in the same stabilizer $G_x$ are at distance one, and two points in distinct stabilizers $G_a$ and $G_b$ (for $a,b\in X^0$) are at distance $d(a,b)$. Let $\pi:Y\to X$ denote the canonical projection.  For two points $x,y\in Y$,  the sets $C_Y(x,y)$ are given by $\pi^{-1}\left(C\left(\pi(x),\pi(y)\right)\right)$ (where $C$ is defined on $X^0$ as in the beginning of the proof). Condition (i) is obviously satisfied, for (ii) the polynomial is $cP(n)$, where $c$ is the uniform bound on the cardinality of the stabilizers and (iii) is satisfied as well because the diameter of $C_Y(x,y)$ differs by at most one from that of $C(\pi(x),\pi(y))$. Intuitively, we blow up the stabilizers to get a free action, and the uniform bound on stabilizers allows for this to be done through a quasi-isometry of multiplicative constant 1, which does not affect conditions (i), (ii) and (iii) of Proposition \ref{max}.\end{proof}
In particular, since Coxeter groups are known to act properly discontinuously by isometries on CAT(0) cube complexes of finite dimension \cite{NibloReeves}, we have the following (already proved in \cite{these} after the suggestion of N. Higson).  
\begin{cor}Coxeter groups have property RD.\end{cor}
Corollary \ref{NoAmSup} is a straight consequence of Theorem \ref{RDcube} and Theorem \ref{obstruction}. Note that even though D. Wise and M. Sageev have a stronger result in \cite{WS}, the techniques are completely different. It was known previously  that virtually soluble subgroups of such groups have to be virtually abelian (see II.7.8 in \cite{BH}). 
\begin{rem}\label{Shahar}The following example has been provided by S. Mozes and shows that the assumption regarding uniform bound on stabilizers cannot be removed. Let $p$ be a prime and ${\bf F}_p$ be a finite field of cardinality $p$. The group $\Gamma=PGL_2({\bf F}_p[t,t^{-1}])$ (the quotient, by the center, of invertible 2 by 2 matrices with coefficients Laurent polynomials in one variable on the finite field ${\bf F}_p$) is generated by the elements 
$$\left(\begin{array}{cc}t&0\\ 0&1\end{array}\right),\left(\begin{array}{cc}1&1\\ 0&1\end{array}\right)\hbox{ and }\left(\begin{array}{cc}1&0\\ 1&1\end{array}\right).$$ 
Consider the group $G=PGL_2(F_p((t)))\times PGL_2(F_p((t^{-1})))$ with its associated affine Bruhat-Tits' building $X$, a product of two $(p+1)$-regular trees.  The group $\Gamma$ acts properly on the vertices of $X$ via the diagonal embedding of $\Gamma$ into $G$, under which $\Gamma$ is an irreducible lattice.  The stabilizers of this action are the finite subgroups 
$$L_n=\{\left(\begin{array}{cc}1&P(t)\\ 0&1\end{array}\right)| P(t)\hbox{ is a polynomial of degree at most }n\}$$
which are of cardinality $|L_n|=p^{n+1}$. Now, since the element 
$$\left(\begin{array}{cc}1&t^n\\ 0&1\end{array}\right)=\left(\begin{array}{cc}t^n&0\\ 0&1\end{array}\right)\left(\begin{array}{cc}1&1\\ 0&1\end{array}\right)\left(\begin{array}{cc}t^{-n}&0\\ 0&1\end{array}\right)$$ 
has length $2n+1$, there is a positive constant $C$ such that for each $n$, $L_n\subset B(Cn)$, and $G$ cannot have RD because the elements $\chi_n$ given by the characteristic functions of the subgroups $L_n$ have operator norm as follows
$$\|\chi_n\|_{op}\geq\frac{\|\chi_n*\chi_n\|_2}{\|\chi_n\|_2}=\frac{\||L_n|\chi_n\|_2}{\sqrt{|L_n|}}=\sqrt{p^{n+1}}\|\chi_n\|_2,$$
which contradicts inequality (2) of Proposition \ref{RD}.
\end{rem}

\medbreak

\noindent
{\it Acknowledgments.} We thank Shahar Mozes for Remark \ref{Shahar}, pointing out a missing assumption in an earlier false formulation of Theorem \ref{RDcube} and Alain Valette for a thorough read-through of the paper with many helpful comments and suggestions. We would also like to thank Daniel Groves for pointing out a missing step in the proof of Lemma \ref{AlmostL}.  The first author would like to thank Daan Krammer for the winter 2001 discussions on property RD, as they inspired the formulation of Proposition \ref{max}. Part of this article was written up during her visit to TIFR Bombay in Summer 2003, and she is grateful to the institute for offering ideal working conditions.  The second author would like to thank the Tufts University FRAC committee for summer funding in 2003 as well as Guoliang Yu for his comments.  Both authors would like to thank the referees for their helpful comments and suggestions. 

\end{document}